\documentclass[letterpaper, 10 pt, conference]{ieeeconf}  

\IEEEoverridecommandlockouts                              

\usepackage{lipsum}
\usepackage[noadjust]{cite}  
\usepackage{amsmath,amssymb}
\usepackage{aligned-overset}
\usepackage{enumerate}
\usepackage{bm}
\usepackage{amsfonts}
\usepackage{anyfontsize}
\usepackage{cancel}

\usepackage{float}

\usepackage{empheq}

\usepackage{multirow}
\usepackage{array}
\usepackage{tabularx}
\usepackage{booktabs}

\newtheorem{theorem}{Theorem}[section]
\newtheorem{lemma}{Lemma}[section]

\newtheorem{remark}{Remark}[section]

\newtheorem{problem}{Problem}[section]

\newcommand\numberthis{\addtocounter{equation}{1}\tag{\theequation}}
\makeatletter

\usepackage{algpseudocode}
\usepackage{algorithm}

\usepackage{setspace}
    
\usepackage{bm} 

\allowdisplaybreaks  

\usepackage[unicode]{hyperref}
\hypersetup{
    colorlinks,
    citecolor=black,
    filecolor=black,
    linkcolor=black,
    urlcolor=black
}

\usepackage{nccmath} 

\usepackage{graphicx}

\usepackage[font=sf]{caption}
\usepackage{subcaption}
\usepackage{mathrsfs}   

\newcommand{\me}{\mathrm{e}}
\newcommand{\md}{\mathrm{d}}

\usepackage{accents}

\usepackage{nccmath} 

\usepackage{pgfplots}
\pgfplotsset{compat=1.15}
\usetikzlibrary{patterns}
\usetikzlibrary{bending}
\usetikzlibrary{calc}
\usetikzlibrary{arrows,backgrounds,positioning,fit,petri}
\usetikzlibrary{decorations.pathreplacing,decorations.pathmorphing,decorations.markings,decorations.shapes}
\usetikzlibrary{arrows.meta}
\usetikzlibrary{quotes,angles}
\usetikzlibrary{chains}
\usetikzlibrary{hobby}
\usepackage{tkz-euclide}
\definecolor{ttttff}{rgb}{0.2,0.2,1}
\definecolor{ududff}{rgb}{0.30196078431372547,0.30196078431372547,1}
\definecolor{xfqqff}{rgb}{0.4980392156862745,0,1}
\definecolor{qqzzff}{rgb}{0,0.6,1}
\definecolor{ffxfqq}{rgb}{1,0.4980392156862745,0}
\definecolor{qqqqff}{rgb}{0,0,1}
\definecolor{ffqqff}{rgb}{1,0,1}
\definecolor{qqzzqq}{rgb}{0,0.6,0}
\definecolor{ffqqqq}{rgb}{1,0,0}
\definecolor{qqwuqq}{rgb}{0,0.39215686274509803,0}
\definecolor{zzttqq}{rgb}{0.6,0.2,0}
\definecolor{wqwqwq}{rgb}{0.3764705882352941,0.3764705882352941,0.3764705882352941}
\definecolor{qqttcc}{rgb}{0,0.2,0.8}
\definecolor{windBlue}{rgb}{0.3176470588,0.9254901961,0.9450980392}

\newcommand\norm[1]{\| #1 \|}

\renewcommand{\thetable}{\arabic{table}}

\makeatletter
\newcommand{\pushright}[1]{\ifmeasuring@#1\else\omit\hfill$\displaystyle#1$\fi\ignorespaces}
\newcommand{\pushleft}[1]{\ifmeasuring@#1\else\omit$\displaystyle#1$\hfill\fi\ignorespaces}
\makeatother

\title{\LARGE \bf A Bearing-rate-free Implementation of a Scalar Estimator for Bearing-only Target Localization and Circumnavigation}

\author{Donglin Sui$^{1}$ and Mohammad Deghat$^{1}$
\thanks{$^{1}$The authors are with the School of Mechanical and Manufacturing Engineering, University of New South Wales, Sydney NSW 2052, Australia (e-mail: {\tt\small d.sui@unsw.edu.au, m.deghat@unsw.edu.au}).}%
}

\begin{document}

\maketitle
\thispagestyle{empty}
\pagestyle{empty}

\begin{abstract}
    This paper studies the problem of guiding an agent to localize and circumnavigate a stationary target at a prescribed radius using bearing measurements. Direct implementation of a classical scalar signed-distance estimator requires explicit evaluation of bearing rate signals that can be difficult to obtain in practice. To eliminate this dependence, we propose a bearing-rate-free (BRF) realization of the scalar estimator that recovers the signed-distance estimate algebraically, without numerical differentiation or a separate filter for estimating the bearing rate. We prove, under consistent initialization, that the proposed BRF realization inherits the guarantees previously established for the classical realization, including exponential convergence of both the target estimation error and radial tracking error, as well as the agent-target collision avoidance. Simulations show that, with noise-free bearings, the BRF realization achieves numerical accuracy comparable to that of an exact-bearing-rate benchmark and, at the higher noise levels tested, exhibits smaller radial tracking oscillations than the backward difference implementation of the classical estimator.
\end{abstract}


\section{INTRODUCTION}
\label{sec:Introduction}

The need to observe a target from different viewpoints has motivated extensive research on target circumnavigation, which typically requires one or more agents to move around the target while maintaining a prescribed distance from it. Early control designs often assume that the agents have access to the target state. When the target is uncooperative, however, even its position may have to be inferred from onboard measurements. The agents' motion then determines both how closely they follow the desired circumnavigation orbit and what information they acquire about the target. This coupling between navigation and localization gives rise to a challenging dual control problem~\cite{feldbaum1963Dual,rantzer2023Dual}.

The choice of sensor is particularly important when an agent must localize the target without revealing its own presence. An active sensor that measures both bearing and distance to the target typically transmits a probing signal that the target could detect. A passive optical sensor, by contrast, can determine the direction to a visible target from the light it receives. For example, an agent can obtain the bearing to a target from a single monocular image if it can identify the target's boundary. The ability to acquire bearings passively and at low cost has helped motivate research on \emph{bearing-only target localization and circumnavigation} (BoTLC).

\begin{figure}[h]
    \centering
    \scalebox{0.74}{%
        \begin{tikzpicture}[scale=1.5, line cap=round,line join=round,>=triangle 45,x=1cm,y=1cm]
            \clip(-4,-3.4) rectangle (1.2,1.56);
            %
            %
            \coordinate (Y) at (-2.9019083532629173,-2.0826434806438683);
            \coordinate (X) at (0,0);
            \coordinate (Xhat) at ($(Y)!0.5!(X)$);
            \definecolor{seqOrange}{RGB}{255,127,0}
            %
            %
            \coordinate (Yrho) at ($(Y)!1.80cm!90:(X)$);
            \coordinate (Xrho) at ($(X)+(Yrho)-(Y)$);
            \coordinate (Ysplit) at ($(Y)!1.25cm!90:(X)$);
            \coordinate (Xsplit) at ($(X)+(Ysplit)-(Y)$);
            \coordinate (Xhatsplit) at ($(Xhat)+(Ysplit)-(Y)$);
            \coordinate (DimensionTick) at ($(Yrho)!0.045cm!90:(Xrho)-(Yrho)$);
            %
            %
            \draw [line width=0.3pt,color=gray!70] (Y)-- (Yrho);
            \draw [line width=0.3pt,color=gray!70] (X)-- (Xrho);
            \draw [line width=0.3pt,color=gray!70] (Xhat)-- (Xhatsplit);
            %
            %
            \draw [line width=0.4pt] (Yrho)--
                node[midway,fill=white,inner sep=2pt] {\Large $\rho(t)$} (Xrho);
            \foreach \P in {Yrho,Xrho}{
                \draw [line width=0.4pt]
                    ($(\P)+(DimensionTick)$)-- ($(\P)-(DimensionTick)$);
            }
            %
            %
            \draw [line width=0.4pt] (Ysplit)--
                node[midway,fill=white,inner sep=2pt] {\Large $\hat{\rho}(t)$} (Xhatsplit);
            \draw [line width=0.4pt] (Xhatsplit)--
                node[midway,fill=white,inner sep=2pt] {\Large $|\tilde{\rho}(t)|$} (Xsplit);
            \foreach \P in {Ysplit,Xhatsplit,Xsplit}{
                \draw [line width=0.4pt]
                    ($(\P)+(DimensionTick)$)-- ($(\P)-(DimensionTick)$);
            }
            %
            %
            \draw [line width=1pt,dash pattern=on 5pt off 5pt] (Y)-- (X);
            %
            %
            \draw [-{Triangle[length=5mm,width=2mm]},line width=2pt,color=ffqqqq]
                (Y)-- (-2.0894814568450726,-1.499580415464064);
            \draw [color=ffqqqq] (-2.7,-1.05) node[anchor=north west] {\Large $\bm{\varphi}(t)$};
            %
            %
            \draw [-{Triangle[length=5mm,width=2mm]},line width=2pt,color=ffqqqq]
                (Y)-- (-2.318845288083113,-2.895070377061713);
            \draw [color=ffqqqq] (-2.4,-2.8) node[anchor=north west] {\Large $\bar{\bm{\varphi}}(t)$};
            %
            %
            \draw [line width=0.8pt,dotted] (Y)-- (-1.74573,-2.08264);
            %
            \draw [
                shift={(Y)},
                -{Stealth[bend, length=3mm, width=2.5mm]}, 
                line width=2pt,
                color=qqttcc
            ] (0.00017248740479756843:0.5202800621839034)
                arc (0.00017248740479756843:33:0.5202800621839034);
            \draw [shift={(Y)},line width=2pt,color=qqttcc] (0,0)-- (0:0.5202800621839034);
            \draw [color=qqttcc] (-2.3353811744404442,-1.65)
                node[anchor=north west] {\Large $\theta(t)$};
            %
            \draw [fill=black] (Y) circle (3pt);
            \draw (-3.5,-2.2) node[anchor=north west] {\Large $\bm{y}(t)$};
            %
            %
            \draw [line width=1pt,dash pattern=on 4pt off 4pt,color=qqzzqq] (X) circle (1cm);
            \draw [line width=1pt,dotted,color=qqzzqq] (X)-- (0.9097994775058879,0.41504808243143776);
            \draw [color=qqzzqq] (0.23133379900014597,0.6)
                node[anchor=north west] {\Large $d^*$};
            %
            \draw [fill=black] (X) circle (3pt);
            \draw (X) node[anchor=north west] {\Large $\bm{x}$};
            %
            %
            \draw [color=seqOrange,fill=seqOrange] (Xhat) circle (3pt);
            \draw [color=seqOrange] ($(Xhat)+(0.04,-0.04)$)
                node[anchor=north west] {\Large $\hat{\bm{x}}(t)$};
        \end{tikzpicture}%
    }
    \caption{\label{fig:notation_conventions}Problem geometry and notation. Other symbols such as $d^*$, $\theta(t)$, and $\tilde{\rho}(t)$ will be explained in subsequent sections.}
\end{figure}
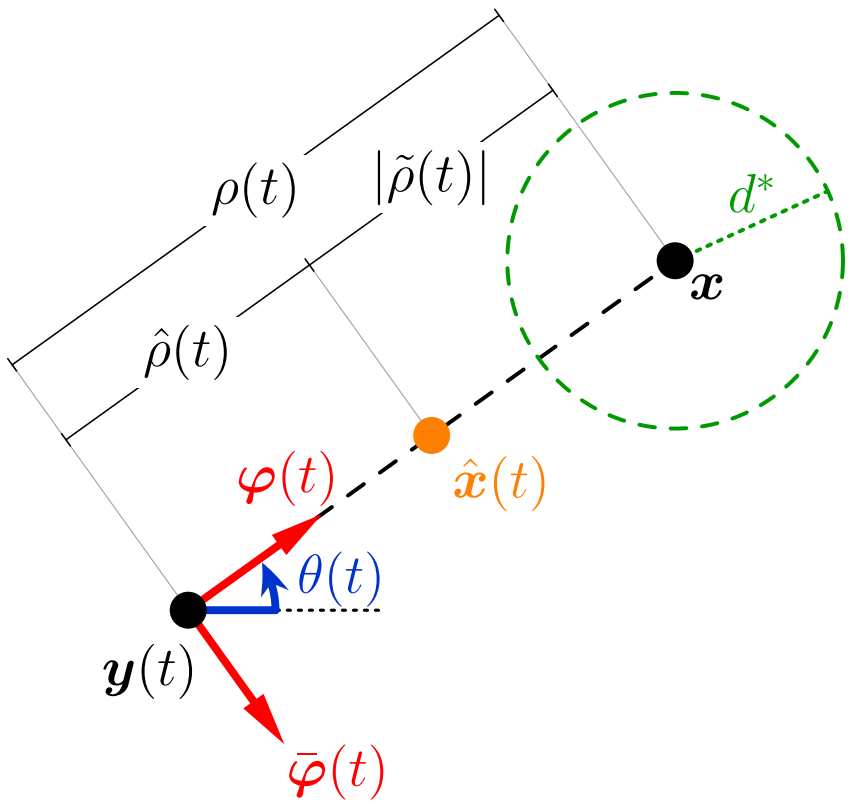

To better explore the current research landscape, we introduce the basic notation. Let $\bm{y}(t)\in\mathbb{R}^2$ and $\bm{x}\in\mathbb{R}^2$ be the positions of the controlled agent and that of the stationary target, respectively. The agent's bearing measurement to the target is represented by the unit vector $\bm{\varphi}(t)$, defined as,
\begin{equation}
    \label{eq:def_varphi}
    \bm{\varphi}(t)=\frac{\bm{x}-\bm{y}(t)}{\norm{\bm{x}-\bm{y}(t)}} =: \frac{\bm{x}-\bm{y}(t)}{\rho(t)} ,
\end{equation}
where $\rho(t)$ denotes the Euclidean distance between the agent and the target. Rotating $\bm{\varphi}(t)$ clockwise by $\pi/2$ gives the perpendicular unit vector $\bar{\bm{\varphi}}(t)\in\mathbb{R}^2$. A graphical illustration is provided in Fig.~\ref{fig:notation_conventions} to aid the discussion. 

An early BoTLC architecture in \cite{deghat2010Target,deghat2014Localization} couples a circumnavigation controller with an adaptive target-state estimator. The controller steers the agent around the target at a specified radius, while the estimator updates the Cartesian estimate of the target position, denoted by $\hat{\bm{x}}(t)\in\mathbb{R}^2$ via the update law \cite[Eq.~(4)]{deghat2010Target}:
\begin{equation}
    \label{eq:MD_estimator}
    \dot{\hat{\bm{x}}}(t) = \gamma \bigl( I_2 - \bm{\varphi}(t)\bm{\varphi}^\top\!(t)\bigr)\bigl( \bm{y}(t)-\hat{\bm{x}}(t)\bigr), \quad \gamma>0,
\end{equation}
where $I_2$ denotes the $2$-by-$2$ identity matrix. Extensions of this vector estimator address the localization of multiple targets~\cite{chun2020Multitargets,ji2024Adaptive}, unbiased localization of a constant velocity target~\cite{sui2024Unbiased}, and circumnavigation by multiple agents~\cite{boccia2017Tracking,dou2020Target,sui2024CommunicationFree,sui2024CollisionFree,jiang2026CommunicationFree}.

Note that during transients, the target position estimate $\hat{\bm{x}}$ produced by the vector estimator \eqref{eq:MD_estimator} need not lie on the current line of sight. The measured bearing, however, already restricts the possible target positions to a ray. If the agent is at $\bm{y}$ and observes the unit bearing $\bm{\varphi}$, then every admissible target position has the form $\bm{y}+\rho\bm{\varphi}$ for some distance $\rho>0$. Cao et al. \cite{cao2020LowCost} are among the first to exploit this additional geometric constraint. Their estimator propagates a scalar signed-distance estimate $\hat{\rho}(t)\in\mathbb{R}$ via the update law \cite[Eq.~(1)]{cao2020LowCost},
\begin{subequations}
    \label{eq:scalar_construction}
    \begin{equation}
        \label{eq:original_scalar_estimator}
        \dot{\hat{\rho}}(t) := - \bm{\varphi}^\top\!(t) \dot{\bm{y}}(t) + k_e \bar{\bm{\varphi}}^\top\!(t) \dot{\hat{\bm{x}}}(t), \quad k_e \geq 1,
    \end{equation}  
    and subsequently reconstructs the target position estimate as
    \begin{equation}
        \label{eq:def_x_hat}
        \hat{\bm{x}}(t) := \bm{y}(t) + \hat{\rho}(t)\bm{\varphi}(t).
    \end{equation}
\end{subequations}
The scalar estimator \eqref{eq:scalar_construction} is particularly useful with multiple targets since each additional target requires only one additional signed-distance estimate \cite{cao2021Safe,wang2024Achieving,zhang2024BearingOnly} and collision avoidance follows readily by construction \cite[Remark~5]{cao2020LowCost}. Extensions of this scalar approach accommodate double-integrator agents~\cite{cao2023Bearingonly,cao2025Fully} and nonholonomic agents~\cite{cao2026Bearingonly}. Finite-time circumnavigation is achieved in \cite{shi2025Finitetime} using a sliding-mode variant of the same scalar estimator idea. Further developments include extensions of the scalar estimator approach to inter-agent distance estimation for cooperative target entrapment~\cite{su2023BearingBased} and adaptive target-entrapping control of leader-follower autonomous underwater vehicles (AUVs)~\cite{su2025Adaptive}. The approach has also been extended to admissible non-circular orbits while guaranteeing a prescribed minimum safe distance between the agent and the target~\cite{wang2026Safe}.

However, the classical realization \eqref{eq:original_scalar_estimator} requires explicit evaluation of bearing rate information. This rate dependence arises through the term $k_e\bar{\bm{\varphi}}^\top\!(t)\dot{\hat{\bm{x}}}(t)$ in \eqref{eq:original_scalar_estimator}, since differentiating the Cartesian reconstruction $\hat{\bm{x}}$ in \eqref{eq:def_x_hat} introduces the bearing derivative $\dot{\bm{\varphi}}(t)$. This can be undesired because obtaining bearing rate is a practical difficulty in vision-based implementations. For example, while a single monocular image can readily provide the bearing from an agent to a target, determining how fast that bearing changes is more demanding as differencing directions extracted from successive frames can amplify noise in the measured pixel coordinates \cite{moshtagh2006Visionbased,jana2022Interception}. In such cases, one possible numerical implementation uses backward difference methods to approximate the rate-dependent terms required by the scalar estimator \eqref{eq:original_scalar_estimator}. However, as noted in \cite{wang2026Safe}, numerical differentiation using finite differences changes the nominal estimation error dynamics, so the continuous-time convergence and safety guarantees do not automatically carry over.

Therefore, the preceding discussion prompts us to revisit the classical scalar estimator of Cao et al.~\cite{cao2020LowCost} with the aim of reproducing its estimates without explicitly evaluating the bearing rate or any rate-dependent terms. The main contributions are summarized as follows.
\begin{itemize}
    \item By introducing a coordinate transformation, we propose a bearing-rate-free (BRF) realization of the classical scalar estimator that recovers the signed-distance estimate $\hat{\rho}$ \emph{algebraically}. In contrast to the existing studies \cite{cao2021Safe,wang2024Achieving,zhang2024BearingOnly,cao2023Bearingonly,cao2025Fully,cao2026Bearingonly,shi2025Finitetime,su2023BearingBased,su2025Adaptive,wang2026Safe}, the proposed realization requires \emph{neither} numerical  differentiation \emph{nor} a memory-based filter for rate estimation.
    
    \item We prove rigorously that, under consistent initialization, the proposed BRF realization produces signed-distance and target position estimates identical to those of the classical rate-dependent realization \eqref{eq:original_scalar_estimator} in \cite{cao2020LowCost}. Consequently, under the hypotheses of the classical result, the proposed BRF realization inherits the convergence and collision avoidance guarantees, including exponential convergence of both the target estimation error and radial tracking error. Unlike \cite{wang2026Safe}, which establishes safety guarantees and ultimate error bounds under additional sufficient conditions and regularity assumptions, our result provides exact continuous-time equivalence without introducing a derivative approximation error.
\end{itemize}

\section{PRELIMINARIES}
\label{sec:Preliminaries}
        
    \subsection{Background and Notation}
    \label{sec:background_and_notations}
        
        Consider a 2D plane containing a stationary target and an agent tasked with estimating the unknown target position from bearing measurements and circumnavigating the target at a prescribed distance $d^*>0$. Following \cite{cao2020LowCost}, we assume the agent obeys single-integrator kinematics, 
        \begin{equation}
            \label{eq:agent_kinematics}
            \dot{\bm{y}}(t) = \bm{u}(t) 
        \end{equation}
        where $\bm{u}(t)$ is the agent's velocity control input. Let the bearing angle $\theta(t)\in\mathbb{R}$ be defined as a continuously unwrapped angle, measured from a fixed reference axis to the unit vector $\bm{\varphi}(t)$, with counterclockwise angles taken as positive and $\theta(0)\in[0,2\pi)$. See Fig.~\ref{fig:notation_conventions} for a graphical illustration of the symbols.

    \subsection{Baseline Estimator and Bearing-Rate Dependence}
    \label{sec:baseline_and_BR_dependence}

        In this section, we demonstrate that the original formulation of the scalar target estimator proposed in \cite{cao2020LowCost}, which is widely adopted in subsequent studies (as in the cohort of \cite{cao2021Safe,wang2024Achieving,zhang2024BearingOnly,cao2023Bearingonly,cao2025Fully,cao2026Bearingonly,shi2025Finitetime,su2023BearingBased,su2025Adaptive,wang2026Safe}), necessitates the evaluation of bearing rate signal for implementation. To this end, we recall from \cite{cao2020LowCost} that the original estimator realization \eqref{eq:original_scalar_estimator} is coupled with the circumnavigation controller given by \cite[Eq.~(2)]{cao2020LowCost},
        \begin{equation}
            \label{eq:def_controller}
            \bm{u}(t) = k_c \bigl( \hat{\rho}(t) - d^*\bigr) \bm{\varphi}(t) + k_\omega \bar{\bm{\varphi}}(t) ,
        \end{equation}
        where $k_c,k_\omega>0$ are control gains. Although written implicitly, we show that the term $k_e\bar{\bm{\varphi}}^\top\dot{\hat{\bm{x}}}(t)$ in \eqref{eq:original_scalar_estimator} necessitates bearing rate evaluation. Following the same ideas as in the proof of \cite[Lemma~4.3]{sui2024Unbiased}, we express $\bm{\varphi}(t)$ in terms of the bearing angle $\theta(t)$, as,
        \begin{equation}
            \label{eq:varphi_in_theta}
            \bm{\varphi}(t) = \begin{bmatrix*}
                \cos \theta(t) \\
                \sin \theta(t)
            \end{bmatrix*},
        \end{equation}
        and since $\bar{\bm{\varphi}}(t)$ is obtained by rotating $\bm{\varphi}(t)$ clockwise by $\pi/2$, we can similarly express it in terms of $\theta(t)$:
        \begin{equation}
            \label{eq:varphi_bar_in_theta}
            \bar{\bm{\varphi}}(t) = \begin{bmatrix*}[r]
                \cos \frac{\pi}{2} & \sin \frac{\pi}{2} \\[4pt]
                -\sin \frac{\pi}{2} & \cos \frac{\pi}{2}
            \end{bmatrix*} \begin{bmatrix*}
                \cos \theta(t) \\
                \sin \theta(t) 
            \end{bmatrix*} = \begin{bmatrix*}[r]
                 \sin \theta(t) \\
                -\cos \theta(t)
            \end{bmatrix*}.
        \end{equation}
        By the chain rule, we obtain
        \begin{align}
            \label{eq:varphi_dot_identity}
            \dot{\bm{\varphi}}(t) \overset{\eqref{eq:varphi_in_theta}}{=} \dot{\theta}(t) \begin{bmatrix*}[r]
                - \sin \theta(t) \\
                  \cos \theta(t)
            \end{bmatrix*} \overset{\eqref{eq:varphi_bar_in_theta}}{=} - \dot{\theta}(t) \bar{\bm{\varphi}}(t).
        \end{align}
        Further, the following identity holds:
        \begin{align*}
            \bar{\bm{\varphi}}^\top\!(t)\dot{\bm{\varphi}}(t) \overset{\eqref{eq:varphi_dot_identity}}&{=} \bar{\bm{\varphi}}^\top\!(t)\left(-\dot{\theta}(t)\bar{\bm{\varphi}}(t)\right) \\
            \overset{\eqref{eq:varphi_bar_in_theta}}&{=} - \dot{\theta}(t) \begin{bmatrix*}
                \sin\theta(t) & -\cos\theta(t)
            \end{bmatrix*} \begin{bmatrix*}[r]
                \sin \theta(t) \\
                -\cos \theta(t)
            \end{bmatrix*} \\
            &=-\dot{\theta}(t) \bigl( \sin^2\theta(t) + \cos^2\theta(t)\bigr) = - \dot\theta(t). \numberthis \label{eq:transpose_identity}
        \end{align*}
        For notational brevity, we omit the time argument hereafter whenever no ambiguity could arise from the context. We now show that the original formulation in \eqref{eq:original_scalar_estimator} necessitates rate information,
        \begin{align*}
            \dot{\hat{\rho}} \overset{\eqref{eq:agent_kinematics},\eqref{eq:original_scalar_estimator}}&{=} - \bm{\varphi}^\top\bm{u} + k_e \bar{\bm{\varphi}}^\top \dot{\hat{\bm{x}}} \\
            \overset{\eqref{eq:def_x_hat}}&{=} - \bm{\varphi}^\top \bm{u} + k_e \bar{\bm{\varphi}}^\top \bigl( \dot{\bm{y}} + \dot{\hat{\rho}}\bm{\varphi}+\hat{\rho}\dot{\bm{\varphi}} \bigr) \\
            \overset{\eqref{eq:agent_kinematics}}&{=} - \bm{\varphi}^\top \bm{u} + k_e \bar{\bm{\varphi}}^\top \bigl( \bm{u} + \dot{\hat{\rho}}\bm{\varphi}+\hat{\rho}\dot{\bm{\varphi}} \bigr) \\
            \overset{\bar{\bm{\varphi}}^\top\!\bm{\varphi}\equiv 0}&{=} - \bm{\varphi}^\top \bm{u} + k_e \bar{\bm{\varphi}}^\top \bm{u} + k_e \hat{\rho} \bar{\bm{\varphi}}^\top \dot{\bm{\varphi}} \\
            \overset{\eqref{eq:transpose_identity}}&{=} - \bm{\varphi}^\top \bm{u} + k_e \bar{\bm{\varphi}}^\top \bm{u} - k_e \hat{\rho}\, \dot{\theta}, \numberthis \label{eq:rate_dependence_steps} \\
            \overset{\eqref{eq:def_controller}}&{=} -k_c \hat{\rho} + k_c d^* + k_e k_\omega - k_e \hat{\rho} \, \dot{\theta}. \numberthis \label{eq:rate_dependency}
        \end{align*}
        Therefore, under the original formulation in \cite{cao2020LowCost}, direct implementation of either \eqref{eq:original_scalar_estimator} or \eqref{eq:rate_dependence_steps} requires explicit evaluation of $\dot{\bm{\varphi}}$ or $\dot{\theta}$, respectively. With this in mind, we are ready to state the control problems.

    \subsection{Problem Formulation}
    \label{sec:problem_statements}
        \begin{problem}[Equivalent BRF Realization]
            \label{problem:Localization_Problem}
            Construct a realization of the scalar estimator \eqref{eq:original_scalar_estimator} that produces the same Cartesian reconstruction $\hat{\bm{x}}=\bm{y}+\hat{\rho}\bm{\varphi}$ without differentiating bearing measurements or explicitly evaluating bearing rates. The BRF realization should, under consistent initialization, produce the same $\hat{\rho}$-trajectory as the original estimator \eqref{eq:original_scalar_estimator}. In addition, when coupled with the circumnavigation controller \eqref{eq:def_controller}, the position and distance estimation errors, defined by
            \begin{align}
                \tilde{\bm{x}}(t)
                    &:= \hat{\bm{x}}(t)-\bm{x},
                    \label{eq:def_x_tilde}\\
                \tilde{\rho}(t)
                    &:= \hat{\rho}(t) - \rho(t),
                    \label{eq:def_rho_tilde}
            \end{align}
            should converge to zero exponentially fast.
        \end{problem}

        \vspace{0.25cm}
        
    \begin{problem}[Target Circumnavigation]
        \label{problem:Circumnavigation_problem}
        Establish that the BRF realization in \emph{Problem~\ref{problem:Localization_Problem}}, coupled with \eqref{eq:def_controller}, preserves the original circumnavigation and collision avoidance guarantees under the same hypotheses, gain selection, and consistent initialization:
        \begin{equation}
            \label{eq:same_cond}
            \rho(0)>0, \quad \hat{\rho}(0)\leq \rho(0), \quad k_e\geq 1, \quad k_c >0, \quad k_\omega >0 .
        \end{equation}
        In addition, the radial tracking error, defined by
        \begin{equation}
            \label{eq:def_delta}
            \delta(t) := \rho(t) - d^*,
        \end{equation}
        should converge to zero exponentially fast as established in \cite[Theorem~2]{cao2020LowCost} and satisfy the safety bound in \cite[Theorem~3]{cao2020LowCost}, that is,
        \begin{equation}
            \label{eq:def_CA_goal}
            \rho(t)\geq\min\{d^*,\rho(0)\}>0,
            \quad \forall\,t\geq0.
        \end{equation}
    \end{problem}

    \vspace{\baselineskip}

\section{PROPOSED REALIZATION}
\label{sec:proposed_algorithm}

    To eliminate the dependence on rate information, we propose an \emph{exactly equivalent} realization that avoids evaluating any bearing rate signal. Introduce the following integrating-factor coordinate,
    \begin{subequations}
        \label{eq:BRF_implementation}
        \begin{align}
            \label{eq:z_dot_def}
            \dot{z} &= \me^{k_e \theta} \left(-\bm{\varphi}^\top\bm{u}+k_e\bar{\bm{\varphi}}^\top \bm{u} \right), \\
            z(0) &:= \me^{k_e \theta(0)}\hat{\rho}(0),
        \end{align}
        whose significance will be explained later in \emph{Lemma~\ref{lem:same_input_output}}.
        Then, the proposed BRF realization is
        \begin{equation}
            \label{eq:hat_rho_BRF}
            \hat{\rho}:= \me^{-k_e\theta}z.
        \end{equation}
    \end{subequations}
    Observe that neither $\dot{\theta}$ nor $\dot{\bm{\varphi}}$ is evaluated in \eqref{eq:BRF_implementation}. Under the circumnavigation controller \eqref{eq:def_controller}, \eqref{eq:z_dot_def} can be written as
    \begin{equation}
        \label{eq:explicit_z_dynamics}
        \dot{z} = -k_c z + \bigl( k_c d^* + k_e k_\omega \bigr) \me^{k_e\theta},
    \end{equation}
    and note that, in contrast to \eqref{eq:rate_dependency}, the $z$-dynamics \eqref{eq:explicit_z_dynamics} does not require any explicit evaluation of rate-dependent signals.

    \vspace{\baselineskip}

\section{Equivalence and Preservation of Convergence}
\label{sec:convergence_analysis}
    We first establish in \emph{Lemma~\ref{lem:same_input_output}} that, under consistent initialization, the proposed BRF realization reproduces the $\hat{\rho}$-trajectories of the original estimator. Then, we show in \emph{Theorem~\ref{thm:stability_holds}} that coupling the two realizations with \eqref{eq:def_controller} leads to the same closed-loop trajectories under matched initial conditions, and it follows that the corresponding convergence results of \cite{cao2020LowCost} are preserved.

    \begin{lemma}
        \label{lem:same_input_output}
        Fix any $T>0$. Consider common signals $\bm{y}$, $\bm{u}=\dot{\bm{y}}$, and $\bm{\varphi}$ on $[0,T]$. The transformation
        \begin{equation}
            \label{eq:def_transformation}
            z(t)=\me^{k_e\theta(t)}\hat{\rho}(t)
        \end{equation}
        establishes a one-to-one correspondence between solutions of the original realization \eqref{eq:original_scalar_estimator} and the auxiliary dynamics \eqref{eq:z_dot_def} on $[0,T]$. Further, under consistent initialization, the reconstruction \eqref{eq:hat_rho_BRF} produces the same signed-distance estimate $\hat{\rho}$ and Cartesian estimate $\hat{\bm{x}}=\bm{y}+\hat{\rho}\bm{\varphi}$ as the original realization \eqref{eq:original_scalar_estimator}, throughout $[0,T]$.
    \end{lemma}

    \begin{proof}
        Define the following scalar auxiliary function,
        \begin{equation}
            \label{eq:def_aux_f}
            f(t):= -\bm{\varphi}^\top\!(t) \bm{u}(t) + k_e \bar{\bm{\varphi}}^\top\!(t)\bm{u}(t).
        \end{equation}
        Then, the original realization \eqref{eq:original_scalar_estimator} can be written in terms of $f$, as follows,
        \begin{equation}
            \label{eq:original_rho_hat_dynamics}
            \dot{\hat{\rho}} \overset{\eqref{eq:rate_dependence_steps},\eqref{eq:def_aux_f}}{=} f - k_e \dot{\theta} \hat{\rho} .
        \end{equation}
        Now consider the BRF realization \eqref{eq:BRF_implementation}, whose time derivative can be obtained as,
        \begin{align*}
            \frac{\md}{\md t}\, \hat{\rho} \overset{\eqref{eq:hat_rho_BRF}}{=} \frac{\md}{\md t}\left(\me^{-k_e\theta }z \right) &= \me^{-k_e\theta}\dot{z} -k_e\dot{\theta} \me^{-k_e \theta} z \\
            \overset{\eqref{eq:z_dot_def},\eqref{eq:def_aux_f}}&{=} f - k_e \dot{\theta} \me^{-k_e\theta} z \\
            \overset{\eqref{eq:hat_rho_BRF}}&{=} f - k_e \dot{\theta} \hat{\rho}, \numberthis \label{eq:BRF_rho_hat_dynamics}
        \end{align*}
        which is identical to \eqref{eq:original_rho_hat_dynamics}. Conversely, differentiating $z=\me^{k_e\theta}\hat{\rho}$ along \eqref{eq:original_rho_hat_dynamics} gives \eqref{eq:z_dot_def}. Further, observe that the proposed transformation \eqref{eq:def_transformation} is invertible at every finite time at which $\theta$ is finite. Indeed, for every finite real-valued $\theta$, the term $\me^{k_e\theta}$ is finite and strictly positive, and hence the transformation \eqref{eq:def_transformation} is invertible. Therefore, under consistent initialization, uniqueness of scalar ordinary differential equations (ODEs) implies equality of the $\hat{\rho}$-trajectories. Equality of the target position estimates then follows directly from \eqref{eq:def_x_hat}.
    \end{proof}

    \begin{remark}
        Note that since the bearing angle $\theta(t)$ is defined as an unwrapped angle,  the proposed BRF realization \eqref{eq:BRF_implementation} is well defined on every finite interval under the stated hypotheses. Its auxiliary state $z(t)$ may, nevertheless, grow with accumulated bearing rotation. Indeed, during sustained circumnavigation, the unwrapped bearing angle $\theta$ increases without bound, and the auxiliary coordinate $z=\me^{k_e\theta}\hat{\rho}$ is unbounded even if $\tilde{\rho}$ is exponentially convergent. This observation motivates a coordinate rescaling scheme for prolonged operation, which is left for future work.
    \end{remark}

    \begin{theorem}
        \label{thm:stability_holds}
        Under \eqref{eq:same_cond} and consistent initialization, the proposed BRF realization \eqref{eq:BRF_implementation}, coupled with \eqref{eq:def_controller}, preserves the guarantees of \cite[Theorems~1--3]{cao2020LowCost}, that is, $\tilde{\bm{x}}(t)$ and $\rho(t)-d^*$ converge to zero exponentially fast, and the safety bound \eqref{eq:def_CA_goal} holds for all $t \geq 0$.
    \end{theorem}

    \begin{proof}
        Fix gains $k_e,k_c,k_\omega$ and an arbitrary initial condition
        \begin{equation}
            \label{eq:initial_state}
            \Gamma_0:=(\bm{y}_0,r_0) :=\bigl(\bm{y}(0),\hat{\rho}(0)\bigr)
        \end{equation}
        satisfying \eqref{eq:same_cond}. Let
        \begin{equation}
            \label{eq:original_CL_solution}
            \Gamma(t,\Gamma_0):=\bigl(\bm{y}(t,\Gamma_0),\hat{\rho}(t,\Gamma_0) \bigr)
        \end{equation}
        denote the closed-loop solution obtained by coupling the agent kinematics \eqref{eq:agent_kinematics} with the original circumnavigation controller \eqref{eq:def_controller} and the original realization of the scalar estimator \eqref{eq:original_scalar_estimator}. With the measured bearing $\bm{\varphi}$ defined by \eqref{eq:def_varphi}, the corresponding initial value problem (IVP) is
        \begin{equation}
            \label{eq:initial_value_problem}
            \left\{
            \begin{aligned}
                \dot{\bm{y}} &= \bm{u},\\
                \bm{u} &= k_c(\hat{\rho}-d^*)\bm{\varphi}
                          + k_\omega\bar{\bm{\varphi}},\\
                \dot{\hat{\rho}} &= -\bm{\varphi}^{\top}\dot{\bm{y}}
                            + k_e\bar{\bm{\varphi}}^{\top}\dot{\hat{\bm{x}}},\\
                \bigl(\bm{y}(0),\hat{\rho}(0)\bigr) &= \Gamma_0.
            \end{aligned}
            \right.
        \end{equation}
        By \cite[Theorems~1--3]{cao2020LowCost}, $\Gamma(\cdot,\Gamma_0)$
        is defined on $[0,\infty)$ and satisfies, for all $t\geq0$,
        \begin{gather}
            \rho(t) \geq \min\{d^*,\rho(0)\}>0,  \label{eq:Gamma_rho_lower_bound} \\
            \|\tilde{\bm{x}}(t)\|+|\delta(t)| \leq C\me^{-\lambda t},
            \label{eq:Gamma_exp_bound}
        \end{gather}
        where the constants $C,\lambda>0$ may depend on $\Gamma_0$.

        We now proceed to show the boundedness of the unwrapped bearing angle $\theta$ along $\Gamma(\cdot,\Gamma_0)$ over every finite interval. Fix any $T>0$, and by \eqref{eq:def_varphi}, we get
        \begin{equation}
            \label{eq:bearing_identity}
            \bm{x}-\bm{y}=\rho\bm{\varphi},
        \end{equation}
        and since the target is stationary, differentiating \eqref{eq:bearing_identity} gives
        \begin{align*}
            \dot{\bm{x}}-\dot{\bm{y}}
            &=
            \dot{\rho}\bm{\varphi}
            +\rho\dot{\bm{\varphi}}
            \\
            &\overset{\eqref{eq:varphi_dot_identity}}{=}
            \dot{\rho}\bm{\varphi}
            -\rho\dot{\theta}\bar{\bm{\varphi}} \\
            \overset{\dot{\bm{x}}=0}{\implies} \quad -\dot{\bm{y}} &= \dot{\rho}\bm{\varphi}
            -\rho\dot{\theta}\bar{\bm{\varphi}} \\
            \overset{\eqref{eq:agent_kinematics}:\;\dot{\bm{y}}=\bm{u}}{\implies} \quad -\bm{u} &= \dot{\rho}\bm{\varphi}-\rho \dot{\theta}\bar{\bm{\varphi}} \\
            \implies \quad -\bar{\bm{\varphi}}^\top\!\bm{u}
            &=
            \dot{\rho}
            \bar{\bm{\varphi}}^\top\!\bm{\varphi}
            -\rho\dot{\theta}
            \bar{\bm{\varphi}}^\top\!\bar{\bm{\varphi}} \\
            &=
            -\rho\dot{\theta}.
        \end{align*}
        Subsequently, we obtain
        \begin{align*}
            \rho\dot{\theta}
            &=
            \bar{\bm{\varphi}}^\top\!\bm{u}
            \\
            \overset{\eqref{eq:def_controller}}&{=}
            k_c(\hat{\rho}-d^*)
            \bar{\bm{\varphi}}^\top\!\bm{\varphi}
            +k_\omega
            \bar{\bm{\varphi}}^\top\!\bar{\bm{\varphi}} \\
            &=k_\omega. \numberthis \label{eq:bearing_angular_rate}
        \end{align*}
        Applying \eqref{eq:Gamma_rho_lower_bound} along $\Gamma(\cdot,\Gamma_0)$ gives
        \begin{equation}
            \label{eq:angular_rate_bound}
            0<\dot{\theta}(t)
            \overset{\eqref{eq:bearing_angular_rate}}{=}\frac{k_\omega}{\rho(t)}
            \overset{\eqref{eq:Gamma_rho_lower_bound}}{\leq}
            \frac{k_\omega}{\min\{d^*,\rho(0)\}}.
        \end{equation}
        Hence, for any $T>0$ and $0\leq t\leq T$, it holds that
        \begin{equation}
            \label{eq:theta_is_finite_on_finite_intervals}
            0\leq\theta(t)-\theta(0)
            =\int_0^t\dot{\theta}(s)\; \md s
            \leq\frac{k_\omega T}{\min\{d^*,\rho(0)\}}.
        \end{equation}
        It then follows immediately from \eqref{eq:theta_is_finite_on_finite_intervals} that $\me^{k_e\theta(t)}$ is also finite and positive on the interval $[0,T]$.

        Now consider the IVP for the same agent and controller with the BRF realization \eqref{eq:BRF_implementation}. Let $\theta_0:=\theta(0)$ denote the initial bearing angle at $\bm{y}_0:=\bm{y}(0)$. With consistent initialization, the IVP of the BRF realization is
        \begin{equation}
            \label{eq:BRF_initial_value_problem}
            \left\{
            \begin{aligned}
                \dot{\bm{y}} &= \bm{u},\\
                \bm{u} &= k_c(\hat{\rho}-d^*)\bm{\varphi}
                          + k_\omega\bar{\bm{\varphi}},\\
                \dot z &= \me^{k_e\theta}
                        \bigl(-\bm{\varphi}^{\top}\bm{u}
                        +k_e\bar{\bm{\varphi}}^{\top}\bm{u}\bigr),\\
                \hat{\rho} &= \me^{-k_e\theta}z,\\
                (\bm{y}(0),z(0)) &= (\bm{y}_0,\me^{k_e\theta_0}r_0),
            \end{aligned}
            \right.
        \end{equation}
        where $\bm{\varphi}$ is obtained by evaluating \eqref{eq:def_varphi} along the position trajectory $\bm{y}$ of this IVP, and $\bar{\bm{\varphi}}$ is obtained accordingly using \eqref{eq:varphi_bar_in_theta}. The corresponding continuously unwrapped bearing angle $\theta$ satisfies $\theta(0)=\theta_0$. Let
        \begin{equation}
            \label{eq:BRF_solution}
            \bar{\Gamma}(t,\Gamma_0) = \bigl(\bm{y}(t),\me^{-k_e\theta(t)}z(t) \bigr)
        \end{equation}
        denote the BRF trajectory expressed in the reconstructed coordinates $(\bm{y},\hat{\rho})$. To prove that $\bar{\Gamma}(t,\Gamma_0)=\Gamma(t,\Gamma_0)$ on $[0,T]$, we construct a candidate trajectory $\Omega(t,\Gamma_0)$ for the IVP \eqref{eq:BRF_initial_value_problem}, and then show using uniqueness of solutions that $\Omega$, expressed in the reconstructed coordinates $(\bm{y},\hat{\rho})$, must be identically equal to the BRF solution $\bar{\Gamma}$. To this end, we first take the position trajectory $\bm{y}(\cdot,\Gamma_0)$ from $\Gamma$ as the candidate agent motion in $\Omega$. Then, at every $t\in[0,T]$, the bearing from this candidate agent position $\bm{y}(\cdot,\Gamma_0)$ to the fixed target position $\bm{x}$ is exactly the bearing trajectory in $\Gamma$, since
        \begin{equation*}
            \frac{\bm{x}-\bm{y}(t,\Gamma_0)}{\|\bm{x}-\bm{y}(t,\Gamma_0)\|} \overset{\eqref{eq:def_varphi}}{=} \bm{\varphi}(t,\Gamma_0), \quad \forall \, t\in [0,T].
        \end{equation*}
        While there can be infinitely many $\theta$-trajectories, differ only by a constant multiple of $2\pi$, that correspond to the bearing trajectory $\bm{\varphi}(t,\Gamma_0)$, under consistent initialization, there exists a unique angle trajectory $\theta(t,\Gamma_0)$ with $\theta(0,\Gamma_0)=\theta_0$ along the candidate motion $\bm{y}(t,\Gamma_0)$. 

        To construct a candidate trajectory for the auxiliary state $z$, define the coordinate map $H(\vartheta,r):=\me^{k_e\vartheta}r$. Applying $H$ to the trajectories $\bigl( \theta(t,\Gamma_0),\hat{\rho}(t,\Gamma_0) \bigr)$ in $\Gamma$ at every $t\in[0,T]$ gives 
        \begin{equation}
            \label{eq:BRF_candidate_state}
            \zeta(t,\Gamma_0):=
            H\bigl(\theta(t,\Gamma_0),\hat{\rho}(t,\Gamma_0)\bigr).
        \end{equation}
        Observe that the candidate trajectory $\zeta$ satisfies 
        \begin{align*}
            \zeta(0,\Gamma_0) &= \me^{k_e\theta_0}r_0, \\
            \me^{-k_e \theta(t,\Gamma_0)} \zeta(t,\Gamma_0) &= \hat{\rho}(t,\Gamma_0).
        \end{align*}
        Now consider a candidate control-input trajectory $\bm{u}(t,\Gamma_0)$, taken directly from $\Gamma$. With $\bm{\varphi}(t,\Gamma_0)$ and $\bm{u}(t,\Gamma_0)$, \emph{Lemma~\ref{lem:same_input_output}} ensures that $\zeta(\cdot,\Gamma_0)$ satisfies the BRF state equation \eqref{eq:z_dot_def} on $[0,T]$. To verify that feeding the candidate trajectory $\zeta$ back into the controller gives this same control input, substitute $\zeta$ into the control law \eqref{eq:def_controller} to obtain 
        \begin{align*}
            &\quad k_c\bigl(
                \me^{-k_e\theta(t,\Gamma_0)}\zeta(t,\Gamma_0)-d^*
            \bigr) 
            \bm{\varphi}(t,\Gamma_0)
            +
            k_\omega\bar{\bm{\varphi}}(t,\Gamma_0) \\
            \overset{\eqref{eq:BRF_candidate_state}}&{=}
            k_c\bigl(\hat{\rho}(t,\Gamma_0)-d^*\bigr)
            \bm{\varphi}(t,\Gamma_0)
            +k_\omega\bar{\bm{\varphi}}(t,\Gamma_0) \\
            \overset{\eqref{eq:def_controller}}&{=}\bm{u}(t,\Gamma_0)
            \overset{\eqref{eq:agent_kinematics}}{=} \dot{\bm{y}}(t,\Gamma_0), \quad \forall\, t\in [0,T].
        \end{align*}
        Thus $\Omega:\, t\mapsto(\bm{y}(t,\Gamma_0),\zeta(t,\Gamma_0))$
        solves the IVP \eqref{eq:BRF_initial_value_problem} on $[0,T]$. In addition, for $\rho>0$, the BRF closed-loop dynamics are locally Lipschitz in the states $(\rho,\theta,z)$, so its solution must be unique for the considered initial state $\Gamma_0$. Hence, by uniqueness, the BRF position and auxiliary state trajectories on $[0,T]$ must be $\bm{y}(t,\Gamma_0)$ and $\zeta(t,\Gamma_0)$, respectively. The common position trajectory together with the initial angle $\theta_0$ also determines the same bearing angle trajectory $\theta(t,\Gamma_0)$. It then follows from \eqref{eq:BRF_solution} and \eqref{eq:BRF_candidate_state} that, for every $t\in[0,T]$,
        \begin{align*}
            \bar{\Gamma}(t,\Gamma_0)
            &= \bigl(\bm{y}(t,\Gamma_0),
               \me^{-k_e\theta(t,\Gamma_0)}\zeta(t,\Gamma_0)\bigr) \\
            &= \bigl(\bm{y}(t,\Gamma_0),\hat{\rho}(t,\Gamma_0)\bigr) \\
            &= \Gamma(t,\Gamma_0).
        \end{align*}
        Lastly, since $T$ is arbitrary, the original and BRF realizations generate the same agent position and signed-distance estimate for all $t \geq 0$. Their $\hat{\bm{x}}$-trajectories also coincide everywhere on $[0,T]$ by \eqref{eq:def_x_hat}. The claimed bounds \eqref{eq:Gamma_rho_lower_bound} and \eqref{eq:Gamma_exp_bound} therefore hold for the BRF solution $\bar{\Gamma}(t,\Gamma_0)$ with the same constants $C$ and $\lambda$ as for $\Gamma(\cdot,\Gamma_0)$. Since $\Gamma_0$ is arbitrary, this conclusion holds for every admissible initial state under
        consistent initialization.
    \end{proof}

    \begin{table*}[t]
        \centering
        \caption{Configurations used in \emph{Study 1}. All configurations are coupled to the original circumnavigation controller \eqref{eq:def_controller} with identical gains and initial conditions. BD denotes backward difference; RK4 denotes the fourth-order Runge-Kutta method; FE denotes forward Euler method. At $t_n:=n\Delta t$, the BD approximation is $\bigl(\mathrm{D}_{\Delta t}^{-}\bm{\varphi}\bigr)_n=(\bm{\varphi}_n-\bm{\varphi}_{n-1})/\Delta t$ for $n\geq1$, initialized to $\bm{0}_2$ at $n=0$ and held constant over $[t_n,t_{n+1})$, where $\rm{D}_{\Delta t}^{-}$ is the backward finite-difference operator and $\bm{\varphi}_n:=\bm{\varphi}(n\Delta t)$. The exact-rate configuration and the ground-truth reference are simulation benchmarks only.}
        \label{tab:accuracy_configurations}
        \begingroup
            \small
            \setlength{\tabcolsep}{4pt}
            \renewcommand{\arraystretch}{1.15}
            \begin{tabularx}{\textwidth}{@{}
                >{\raggedright\arraybackslash}p{0.235\textwidth}
                >{\raggedright\arraybackslash}p{0.205\textwidth}
                >{\raggedright\arraybackslash}X
                >{\raggedright\arraybackslash}p{0.105\textwidth}@{}}
            \toprule
            \textbf{Configuration} & \textbf{Estimator/model}
            & \textbf{Bearing-rate treatment} & \textbf{Integrator} \\
            \midrule
            Cao et al.\ (2020) + BD + RK4
            & Original realization \eqref{eq:original_scalar_estimator}
            & BD, held over each sampling interval
            & RK4 \\
            \addlinespace[2pt]
            Cao et al.\ (2020) + BD + FE
            & Original realization \eqref{eq:original_scalar_estimator}
            & BD, held over each sampling interval
            & FE \\
            \addlinespace[2pt]
            BRF + RK4
            & Proposed realization \eqref{eq:BRF_implementation}
            & Not required
            & RK4 \\
            \addlinespace[2pt]
            Cao et al.\ (2020) + exact rate
            & Original realization \eqref{eq:original_scalar_estimator}
            & Obtained using $\dot{\varphi}=-k_\omega\bar{\varphi}/\rho$, evaluated at every RK stage
            using the simulated true agent-target distance $\rho$
            & RK4 \\
            \addlinespace[2pt]
            Ground-truth reference
            & Closed-loop ODEs \eqref{eq:groundtruth}
            & Rate dependence incorporated analytically; no numerical differentiation
            & \texttt{MATLAB} \texttt{ode113} \\
            \bottomrule
            \end{tabularx}
            \par\smallskip
        \endgroup
    \end{table*}

    \begin{figure*}[h]
        \centering
        \includegraphics[width=1.0\linewidth]{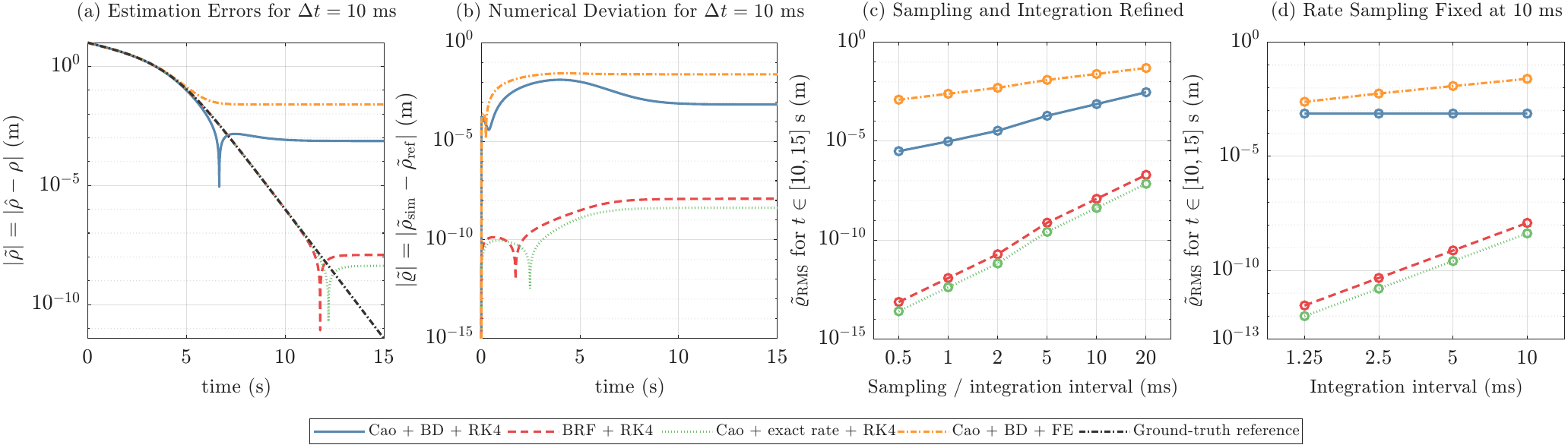}
        \caption{\label{fig:num_acc}Numerical accuracy comparison for \emph{Study 1}.}
    \end{figure*}

    \vspace{\baselineskip}

    \begin{figure*}[h]
        \centering
        \includegraphics[width=1.0\linewidth]{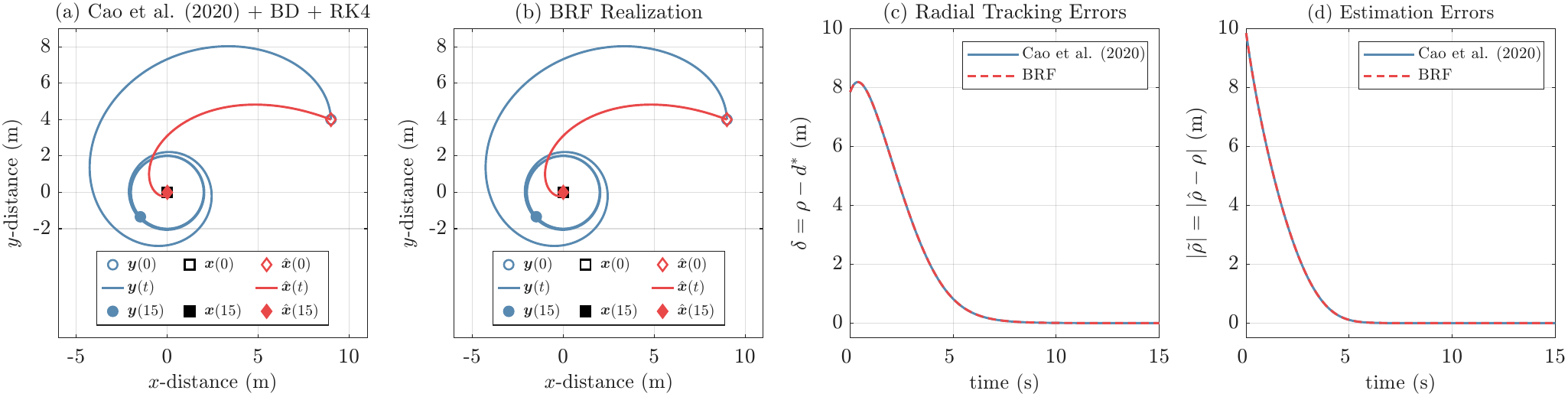}
        \caption{\label{fig:stationary_target_ideal}Simulation plots for \emph{Study 1}. Trajectories are generated under the Cao et al. (2020) + BD + RK4 configuration and the proposed BRF configuration, evaluated at $\Delta t=h=1\, \rm{ms}$.}
    \end{figure*}

    \begin{figure*}[h]
        \centering
        \includegraphics[width=1.0\linewidth]{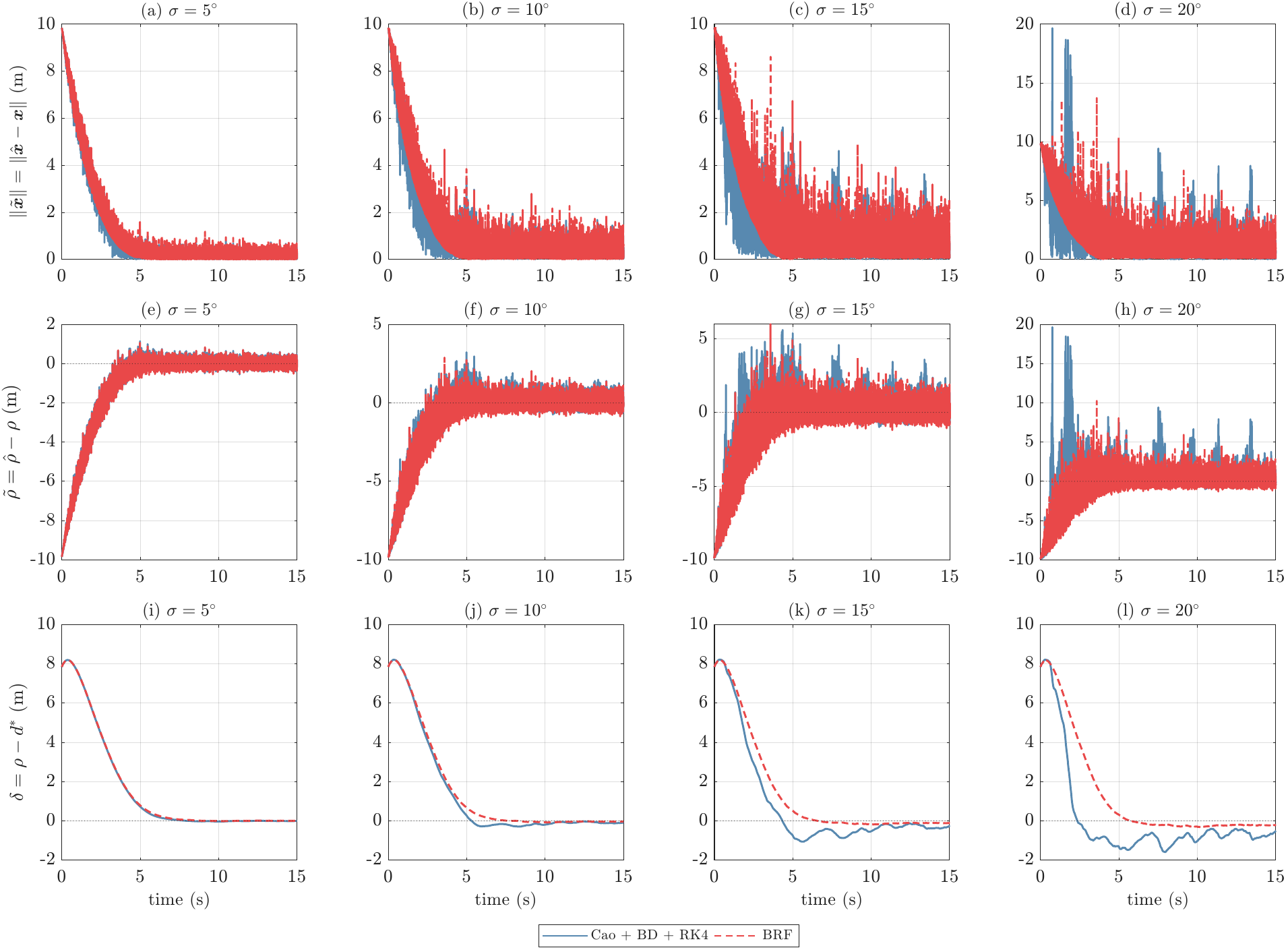}
        \caption{\label{fig:noise_acc_traj}Simulation results for \emph{Study 2} with an RK4 step of $1\,\mathrm{ms}$.}
    \end{figure*}


    \begin{figure}[h]
        \centering
        \includegraphics[width=1.0\linewidth]{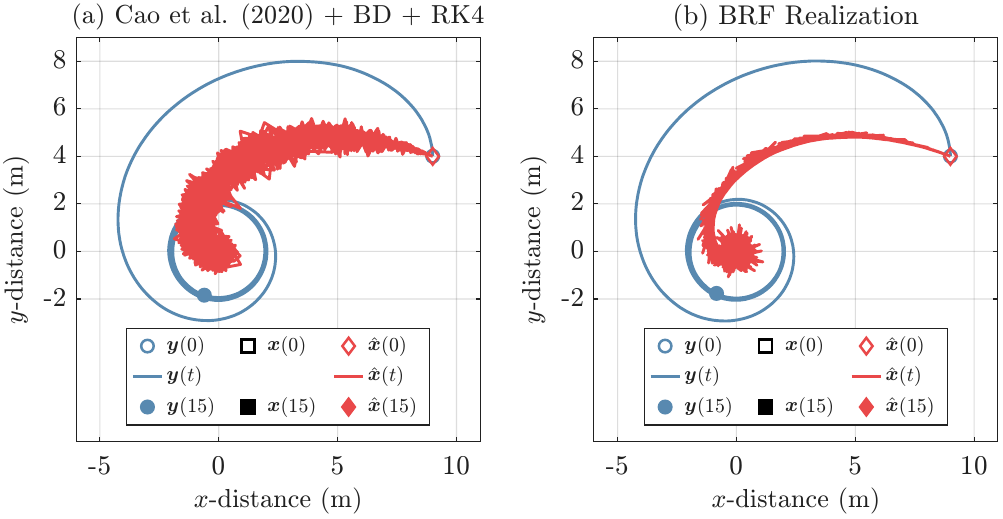}
        \caption{\label{fig:sample_noise_traj}Illustrative agent and target estimate trajectories for $\sigma=5^\circ$ and $\Delta t = 1\,\rm{ms}$.}
    \end{figure}

\section{SIMULATION STUDIES}
\label{sec:Simulations}

    We conduct two simulation studies to compare the numerical accuracy of the original and BRF realizations and to examine their behavior under noise-corrupted bearing measurements. In both studies, we adopt the same initialization and gain selection as in \cite{cao2020LowCost}: the target is stationary at $\bm{x}=[0,0]^\top\,\mathrm{m}$, the agent is initially positioned at $\bm{y}(0)=[9,4]^\top\,\mathrm{m}$ with the agent's estimate initialized to $\hat{\rho}(0)=0$, and the desired circumnavigation radius is set to $d^*=2\,\mathrm{m}$. Throughout the two simulation studies, we use the circumnavigation controller \eqref{eq:def_controller} with $k_c=k_e=1$ and $k_\omega=5$. 

    \subsection{Study 1: Numerical accuracy with noise-free bearings}

    In this simulation study, we compare the four configurations, as summarized in Table~\ref{tab:accuracy_configurations}, when noise-free bearing measurements are available. The original estimator realization in \eqref{eq:original_scalar_estimator} is implemented using a backward difference (BD) approximation of the bearing rate $\dot{\bm{\varphi}}$ and either the fourth-order Runge--Kutta method (RK4) or forward Euler (FE) method. Let $\Delta t$ denote the interval at which the BD approximation is updated and let $h$ denote the integration step. The BD rate is held constant over each update interval, including all integration substeps within that interval.  The proposed BRF realization \eqref{eq:BRF_implementation} is integrated using RK4. An additional implementation of the original estimator realization \eqref{eq:original_scalar_estimator} uses the exact bearing rate at every RK4 stage, which serves as a benchmark representing the case where rate information is directly available.

    To quantify numerical accuracy, we construct a high-accuracy/ground-truth reference by integrating the reduced closed-loop equations
    \begin{subequations}
        \label{eq:groundtruth}
        \begin{empheq}[left=\empheqlbrace]{align}
            \dot{\rho}(t) &= k_c\bigl(d^*-\rho(t)+q(0)\me^{-k_e\vartheta(t)}\bigr), \label{eq:GT_rho_dot} \\
            \dot{\vartheta}(t) &= \frac{k_\omega}{\rho(t)}, \label{eq:GT_vartheta_dot} \\
            q(t) &= q(0)\me^{-k_e\vartheta(t)}, \label{eq:GT_q}
        \end{empheq}
    \end{subequations}
    where $q(t):=-\tilde{\rho}(t)$ and $\vartheta(t):=\theta(t)-\theta(0)$, with $\rho(0)=q(0)=\sqrt{97}\,\mathrm{m}$ and $\vartheta(0)=0$. Equations \eqref{eq:GT_rho_dot}--\eqref{eq:GT_vartheta_dot} are solved using MATLAB's \texttt{ode113}, with relative and absolute tolerances of $3\times10^{-14}$ and $3\times10^{-15}$, respectively, and a maximum step of $0.025\,\mathrm{s}$. Equation \eqref{eq:GT_q} is used to reconstruct the reference distance estimation error as $\tilde{\rho}_{\mathrm{ref}}(t)=-q(t)$. Further, since the estimation error remains nonzero at finite times, its magnitude alone does not measure numerical accuracy. Hence, we define the numerical deviation,
    \begin{equation}
        \label{eq:simulation_numerical_deviation}
        \tilde{\varrho}(t)
        :=\tilde{\rho}_{\mathrm{sim}}(t)-\tilde{\rho}_{\mathrm{ref}}(t),
    \end{equation}
    where $\tilde{\rho}_{\mathrm{sim}}$ denotes the $\tilde{\rho}$-trajectory obtained with a given configuration, whereas $\tilde{\rho}_{\mathrm{ref}}$ denotes the corresponding ground-truth reference trajectory.

    The evolution of the distance estimation error $\tilde{\rho}$ and the numerical deviation $\tilde{\varrho}$ over $[0,15]\,\mathrm{s}$, with a simulation sampling interval of $\Delta t=h=10\,\mathrm{ms}$, is shown in Fig.~\ref{fig:num_acc}(a) and Fig.~\ref{fig:num_acc}(b), respectively. It can be seen that the proposed BRF realization closely follows the reference and achieves numerical accuracy comparable to that of the exact-rate benchmark, whereas both configurations using BD exhibit larger numerical deviations. 

    Subsequently, we distinguish the effects of bearing-rate approximation and numerical integration through two refinement simulations. In Fig.~\ref{fig:num_acc}(c), both simulation intervals are reduced together, with $\Delta t=h\in\{20,10,5,2,1,0.5\}\,\mathrm{ms}$. All configurations become more accurate, and the BRF and exact-rate configurations maintain substantially smaller numerical deviations than the two BD configurations. In Fig.~\ref{fig:num_acc}(d), the BD update interval is fixed at $\Delta t=10\,\mathrm{ms}$, whereas the integration step is reduced to $h\in\{10,5,2.5,1.25\}\,\mathrm{ms}$. The numerical deviation of Cao + BD + RK4 remains nearly unchanged, indicating that its error is dominated by the bearing rate approximation. The FE configuration improves as $h$ decreases, whereas the BRF and exact-rate RK4 configurations continue to approach the reference solution. 

    To provide a more intuitive illustration, Fig.~\ref{fig:stationary_target_ideal} shows a spatial illustration for Cao + BD + RK4 and BRF + RK4 at $\Delta t=h=1\,\mathrm{ms}$. Both agents approach the desired orbit, and their trajectories are nearly indistinguishable, which is consistent with the continuous-time equivalence established in \emph{Lemma~\ref{lem:same_input_output}}.

    \subsection{Study 2: Performance under noise-corrupted bearings}
    
    In this simulation study, we retain the same initialization and setup as in \emph{Study 1} and compare Cao + BD + RK4 with BRF + RK4 under noise-corrupted bearings. Specifically, at each sampling instant $t_n=n\Delta t$, an independent and identically distributed (i.i.d.) Gaussian angular perturbation $\xi_n$ is applied to the true bearing vector $\bm{\varphi}(t_n)$:
    \begin{equation}
        \label{eq:simulation_noise}
        \bm{\varphi}^{\mathrm{m}}_n
        = R(\xi_n)\bm{\varphi}(t_n),
        \quad \xi_n \overset{\text{i.i.d.}}{\sim} \mathcal{N}(0,\sigma^2),
    \end{equation}
    where $R(\cdot)\in\rm{SO}(2)$ is the 2D counterclockwise rotation matrix. We compare the Cao + BD + RK4 and BRF + RK4 configurations over $[0,15]\,\mathrm{s}$ for angular standard deviations $\sigma\in\{5^\circ,10^\circ,15^\circ,20^\circ\}$, with angles converted to radians in the implementation. The noise update interval and the interval used for backward differentiation are fixed at $\Delta t=1\,\mathrm{ms}$, and the same noise-corrupted bearing sequence is used for both configurations. 

    Figure~\ref{fig:noise_acc_traj} compares the target position estimation error $\|\tilde{\bm{x}}\|$, the signed-distance estimation error $\tilde{\rho}$, and the signed radial tracking error $\delta=\rho-d^*$. Note that since $\hat{\bm{x}}$ is reconstructed using the noisy bearing vector as per \eqref{eq:def_x_hat}, the position estimation error $\tilde{\bm{x}}$ includes both distance and angular errors and generally differs from $|\tilde{\rho}|$. In Fig.~\ref{fig:noise_acc_traj}, it can be seen across the columns that the steady-state errors generally increase with the noise level. Figures~\ref{fig:noise_acc_traj}(k)--($\ell$) put in evidence that at the noise level of $\sigma=15^\circ$ and $\sigma=20^\circ$, the BRF realization maintains a significantly less oscillatory radial tracking error, compared to that achieved by the Cao + BD +RK4 configuration.

    To provide a spatial illustration, Fig.~\ref{fig:sample_noise_traj} shows the evolution of the agent and target estimate trajectories for the noise level of $\sigma=5^\circ$. While both configurations approach an orbit close to the desired radius, Fig.~\ref{fig:sample_noise_traj}(b) shows that the estimated target positions under the proposed BRF realization experience less fluctuation in response to the bearing noise.

\section{CONCLUDING REMARKS}
\label{sec:conclusion}

    This paper presents a bearing-rate-free realization of the classical scalar signed-distance estimator in \cite{cao2020LowCost} for localization and circumnavigation of a stationary target. Under consistent initialization, the proposed BRF realization reproduces the original target estimates without explicitly differentiating the measured bearing and preserves exponential convergence of the localization and radial tracking errors. At noise levels tested in simulation, the radial tracking error under the BRF realization exhibited smaller oscillations than that under the original realization with numerical differentiation. This observation also motivates a future theoretical robustness analysis to quantify how the bearing noise level affects the proposed BRF realization. 
    
    Future research directions include extending the analysis to moving targets and investigating whether a bearing-rate-free realization can preserve finite-time convergence under sliding-mode control. A coordinate rescaling scheme should also be investigated for prolonged BoTLC operations.


\addcontentsline{toc}{chapter}{Bibliography}
\bibliographystyle{ieeetr}
\bibliography{myRefs}

\end{document}